\documentclass[preprint,12pt]{elsarticle}
\usepackage{float}
\usepackage{amssymb}
\usepackage{amsmath}
\newenvironment{proof}[1][Proof]{{\em {#1}. }}{{\hfill\raisebox{.6ex}{\framebox[7pt][l]{}}\hspace{-3pt}\vspace{12pt} }}
\newtheorem{theorem}{Theorem}

\usepackage{xcolor}

\makeatletter
\def\ps@pprintTitle{%
  \let\@oddhead\@empty
  \let\@evenhead\@empty
  \let\@oddfoot\@empty
  \let\@evenfoot\@empty
}
\makeatother

\begin{document}

\begin{frontmatter}

\title{Hilbert basis in the face-centered cubic grid}

\author[label1]{B\'ela Vizv\'ari\corref{cor1}}
\author[label2]{Gergely Kov\'acs}
\author[label3,label4]{Benedek Nagy}
\author[label5]{Ne\c{s}et Den\.{i}z Turgay}

\cortext[cor1]{Corresponding author}
\address[label1]{Department of Industrial Engineering, Eastern Mediterranean University, Famagusta 99628, North Cyprus, via Mersin-10, Turkey, e-mail: bela.vizvari@emu.edu.tr}
\address[label2]{Edutus University, 2800 Tatab\'anya, Hungary, e-mail: kovacs.gergely@edutus.hu}
\address[label3]{Department of Mathematics, Faculty of Arts and Sciences,
Eastern Mediterranean University, Famagusta 99628, North Cyprus, via Mersin-10, Turkey, e-mail: benedek.nagy@emu.edu.tr}
\address[label4]{Department of Computer Science, Institute of Mathematics and Informatics, Eszterh\'azy K\'aroly Catholic University, 3300 Eger, Hungary}
\address[label5]{Department of Mathematics, Faculty of Arts and Sciences,
Eastern Mediterranean University, Famagusta 99628, North Cyprus, via Mersin-10, Turkey, e-mail: neset.turgay@emu.edu.tr}

\begin{abstract}
The Hilbert basis is fundamental in describing the structure of the integer points of a polyhedral cone. The face-centered cubic grid is one of the densest packing of the 3-dimensional space. The cycles of a grid satisfy the constraint set of a pointed, polyhedral cone which contains only non-negative integer vectors. The Hilbert basis of a grid gives the structure of the basic cycles in the grid.
It is shown in this paper that the basic cycles of the FCC grid belong to 11 types. It is also discussed that how many elements are contained in the individual types. The proofs of the paper use geometric, combinatorial, algebraic, and operations research methods.
\end{abstract}


\end{frontmatter}

\section{Introduction}\label{introduct}
One of the non-traditional grids in 3D is the face-centered cubic (FCC) grid which also has the point lattice property. This grid plays an important role in the field of crystallography.  The FCC lattice structure is important for understanding the arrangement of atoms in many metals such as  aluminum, copper, silver, and gold since these metals adopt the FCC structure.  The FCC grid is also  an important topic in solid state physics \cite{kittel}. The FCC grid and its higher dimensional generalizations also appear in various theoretical and mathematical studies, e.g. in relation to optimization problems, such as finding minimal distances in networks having any of these structures \cite{IWCIA06,IWCIA09,stom,tur23}.  The FCC grid is the densest packing in 3D. This was proven in \cite{Hales} as part of the Kepler conjecture, providing an efficient way to arrange spheres.
The FCC grid has also applications in other disciplines.
The energy absorption properties of FCC lattice structures are investigated in \cite{wang}.
Recently, scientists investigated FCC cubic carbon as a fourth basic carbon allotrope with properties of intrinsic semiconductors and ultra-wide bandgap \cite{Konyashin}.

Hilbert basis plays an important role in various mathematical and applied contexts. It arises in combinatorics its connection with lecture hall partitions \cite{Olsen}. It is used in integer programming \cite{schrijver98}, in Petri nets \cite{Ma},  in combinatorial optimization \cite{sebo} and to determine the integer Carath\'eodory rank of rational cones \cite{Kuhlmann}. It is utilized to give the Schrijver system for the flow cone in series--parallel graphs \cite{Barbato}. Hilbert basis can also be a tool for understanding total dual dyadicness and dyadic generating sets \cite{Abdi}.
Hilbert basis also arises in graph theory. For instance, class of finite graphs whose sets of edge cuts form Hilbert bases were studied in \cite{Goddyn,Laurent}. Properties of minimal Hilbert bases have also attracted interest for study in the literature \cite{Henk,Liu}.

Hilbert bases for a grid gives important information about the grid and the structure of the set of shortest paths in the grid. In fact, Hilbert bases describes the set of all closed directed cycles in the grid.
In \cite{kov5}, the Hilbert basis is introduced in the body-centered cubic (BCC) grid which is another  important non-traditional grid.  A mathematical relation among various grids, including cubic, the 4-dimensional variant of the body-centered cubic, face-centered cubic and diamond cubic grids  is described, e.g., in \cite{grids}. In this paper, we describe elements of the Hilbert basis in the FCC grid, thus giving important information both about the structure of the grid and structure of the sets of shortest paths that are used to compute digital distances \cite{stom,DAMdist}.  Assume that a directed cycle is divided into two directed paths. If the direction in one of the paths  is changed for the opposite, then the two paths go from the same starting point to the same end point.
The path which has longer length never occur in a minimal path as part of the minimal path.



\section{Description of the Face-Centered Cubic Grid}\label{+++}
\begin{figure}[H]
\centering
   \includegraphics[width=0.3\textwidth]{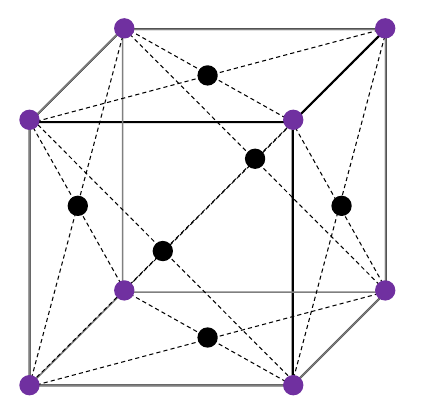}
    \caption{A unit cube of the FCC grid }
 \label{3D-FCC}
\end{figure}

The FCC grid consists of four cubic subgrids: one subgrid contains points with all even-numbered coordinates, while the other three subgrids include points with one even-numbered coordinate and two odd-numbered coordinates. In this grid a point $(x_1,x_2,x_3) \in \mathbb{Z}^3$ has an even number of odd coordinates which tells us that $x_1+x_2+x_3$ is even. From a geometric point of view, the FCC grid consists of unit cubes with a point
at each corner of the cube and a point in the centre of each face of the cube, see Figure \ref{3D-FCC}. This is ``grid space'', and here by a vector we mean a grid vector. From now on when we consider FCC grid we will consider the grid space.

The steps of the FCC grid are represented by the columns of the following matrix.
{\tiny
\begin{center}
\[
\begin{array}{rrrrrrrrrrrrrrrrrrrr}
a_1 & a_2 & a_3 & a_4 & a_5 & a_6 & d_1 & d_2 & d_3 & d_4 & d_5 & d_6 & d_7 & d_8 & d_9 & d_{10} & d_{11} & d_{12} \\
\end{array}
\]
\end{center}
\begin{eqnarray}\label{e001}
\left(
\begin{array}{rrrrrrrrrrrrrrrrrrrr}
 2 & 0 & 0 & -2 & 0 & 0 & 1 & 1 & -1 & -1 & 1 & 1 & -1 & -1 & 0 & 0 & 0 & 0 \\
 0 & 2 & 0 & 0 & -2 & 0 & 1 & -1 & 1 & -1 & 0 & 0 & 0 & 0 & 1 & 1 & -1 & -1 \\
 0 & 0 & 2 & 0 & 0 & -2 & 0 & 0 & 0 & 0 & 1 & -1 & 1 & -1 & 1 & -1 & 1 & -1 \\
\end{array}
\right)
\end{eqnarray}}

The matrix of the FCC grid, the vectors $a_1,...,a_6$, referred to as unit cube vectors, represent steps between two neighbor points of the same cubic grid, and the vectors $d_1,...,d_{12}$, referred to as diagonal vectors,  represent steps between neighbor points of different subgrid. $a$  steps are on the sides of the unit cube, while $d$  steps are diagonal steps between corners and face-centers.  From now on unless otherwise stated,
$A$ denotes the matrix of (\ref{e001}).

\section{Notations}
In what follows, some notations are used. $\mathbb{Z}_{+}$ denotes the set of all non-negative integers. The ordering relation of numbers is expressed by the character $\leq$. The same character is used for the
componentwise ordering of vectors. If $(x_1,x_2,\dots,x_n)^T$ and $(y_1,y_2,\dots,y_n)^T$ are two vectors of real components, then $x\leq y$ means that
$x_i\leq y_i,\,i=1,\dots,n$. Both the zero number and the zero vector are denoted by 0. Letter $e$ denotes a vector such that all of its components are 1.

In what follows, $H$ denotes a Hilbert basis. Hilbert bases are introduced in the next section.

Let $X\subset \mathbb{R}^n$. If $f(x)$ is a real valued function on $\mathbb{R}^n$, then
\begin{eqnarray}\label{lex0}
\min\{f(x)\,:\,x\in X\}
\end{eqnarray}
is an optimization problem which may have more than one optimal solutions. The set of its optimal solutions is denoted by $argmin\{f(x)\,:\,x\in X\}$.
It is worth to mention that {$x^*$ is an optimal solution of (\ref{lex0}) if $x^*\in X$ and $\forall x\in X \,:\, f(x)\geq f(x^*)$.

Multi-objective optimization is a special branch of optimization theory. There are several objective functions on the same set of feasible solutions. The problem is to find a point which has good value in each objective function. Obviously, it is unlikely that a point exists such that it is optimal in all objective function at the same time. There is a need for models that tell us how to consider all objective functions at once. One of the approaches is
the lexicographic ordering of the objective functions. More formally, let $T\subset \mathbb{R}^n$ which is the set of the feasible solutions. Let
$f^T(x)=(f_1(x),f_2(x),\dots, f_m(x))^T$ be a vector valued real function on $T$. The lexicographic approach optimizes the objective functions in index order.
Assume that the direction of the optimization is minimization. The first problem to solve is
\begin{eqnarray}\label{lex1}
\min \left\{ f_1(x):x\in T\right\}.
\end{eqnarray}
Problem (\ref{lex1}) is a normal optimization problem with a single objective function. It can be solved by suitable methods.
Assume that optimal solution exists and the optimal value is $f_1^*$. The first objective function is fixed on this value and the second objective
function is optimized in the second problem which is as follows:
\begin{eqnarray}\label{lex2}
\min \left\{ f_2(x):x\in T,\;f_1(x)=f_1^*\right\}.
\end{eqnarray}
Assume again that an optimal solution exists and the optimal value is $f_2^*$. Then both the first and second objective functions are fixed. The third
problem is
\begin{eqnarray}\label{lex3}
\min \left\{ f_3(x):x\in T,\;f_1(x)=f_1^*,\;f_2(x)=f_2^*\right\}.
\end{eqnarray}
The process is continued in this way. The last problem with single objective function to be solved is
\begin{eqnarray}\label{lexm}
\min \left\{ f_m(x):x\in T,\;f_1(x)=f_1^*,\;f_2(x)=f_2^*,\dots,\;f_{m-1}(x)=f_{m-1}^*\right\}.
\end{eqnarray}
If the lexicographic approach is applied for a multi-objective problem, then the notation of the problem is
\begin{eqnarray}\label{lex}
l-min \left\{ f(x):x\in T\right\}.
\end{eqnarray}

\section{Hilbert Basis in General}\label{HB}

Hilbert introduced an important algebraic notion in \cite{hilbert890}. Later Jeroslow recognized that the same
notion has importance in geometry \cite{jeroslow78}. The treatment below is based on Jeroslow's ideas.

A {\em polyhedral cone} in the  $n$-dimensional Euclidean space is the intersection of finitely many half-spaces.
It means that the polyhedral cone is a set of points $x=(x_1,x_2,\dots,x_n)^T$ satisfying finitely many homogeneous linear inequalities as follows:
\begin{eqnarray}\label{HB01a1}
g_{11}x_1+g_{12}x_2+\cdots +g_{1n}x_n\geq 0\\
\label{HB01a2} g_{21}x_1+g_{22}x_2+\cdots +g_{2n}x_n\geq 0\\ \nonumber\cdots\\
\label{HB01a3} g_{m1}x_1+g_{m2}x_2+\cdots +g_{mn}x_n\geq 0,
\end{eqnarray}
where all the coefficients are real. It is obvious that if the vectors $u$ and $w$ satisfy the inequalities of (\ref{HB01a1})-(\ref{HB01a3}) and $\lambda>0$, then
both $u+w$ and $\lambda u$ satisfy the inequalities of (\ref{HB01a1})-(\ref{HB01a3}) which is the classical definition of convex cone in general. These two properties are called {\em algebraic properties} in this paper. Let

\begin{eqnarray}\label{HB02}
C=\left\{x\in \mathbb{R}^n\mid x\; \text{satisfies} \; (\ref{HB01a1})-(\ref{HB01a3})\right\},
\end{eqnarray}
where all coefficients in (\ref{HB01a1})-(\ref{HB01a3}) are rationals.

The Hilbert basis is a tool to describe the elements of the set
\begin{eqnarray}\label{HB03}
K=C\cap \mathbb{Z}^n.
\end{eqnarray}
The theorem which gives the existence of Hilbert bases is as follows \cite{jeroslow78}:

\begin{theorem}\label{THB1}
If all coefficients in (\ref{HB01a1})-(\ref{HB01a3}) are rationals, then there is a finite subset $H\subset K$ such that every element of
$K$ is a non-negative integer sum of elements of $H$.
\end{theorem}
We will refer these properties as the Hilbert basis property. From now on $C$ is obtained from rational coefficients in (\ref{HB01a1})-(\ref{HB01a3}).
With other words, the theorem claims the existence of a finite subset of $K$ such that
\begin{eqnarray}\label{HB04}
K=\left\{\sum_{h\in H}\lambda_h h\mid \forall h\in H: \lambda_h\in \mathbb{Z}_{+}\right\}.
\end{eqnarray}

The rationality of the coefficients is equivalent to the integrality of the same coefficient. It is a crucial
condition. The construction of \cite{jeroslow71} shows that no such finite subset exists in the irrational case.
The construction can be generalized to a large extent.

If $u\in K\setminus H$ for a Hilbert basis $H$, then $H\cup\{ u\}$ is also a Hilbert basis.

For example the line is a 1-dimensional cone. The set $H_0=\{1,-1\}$ is a Hilbert basis. None of the two elements
can be deleted without losing the Hilbert basis property. However, it is not the only Hilbert basis which is minimal
as a set. Let $p$ and $q$ be two positive and relative prime integers. Then the set $H_{p,q}=\{p,q,-p,-q\}$
is a Hilbert basis such that none of its elements can be deleted without losing the Hilbert basis property.

There are polyhedral cones which have only one minimal Hilbert basis. Convex polyhedrons can be defined in
general as the intersection of finitely many half-spaces. Polyhedrons may have {\em vertex} or in other name {\em extreme point}.
The notion of vertex has equivalent definitions. One of them is that a support plane exists such that it contains
only one point from the polyhedron. Support plane means that the polyhedron is on one side of the plane including the plane itself.
It is obvious that if a cone has a vertex, then it is the origin. The next theorem was discovered
independently both in \cite{jeroslow78} and \cite{schrijver81}.

\begin{theorem}\label{THB2}
If the origin is the vertex of the rational polyhedral cone $C$, then there is a unique minimal Hilbert basis of $C$ such that
every Hilbert basis of $C$ as a set contains completely the elements of the minimal Hilbert basis.
\end{theorem}

The statement of the theorem also means that none of its elements can be deleted without losing the Hilbert basis property.

The only minimal Hilbert basis of the non-negative half-line is {1}. For any finite set $\{g_1,g_2,\dots, g_l \}$ of integers
with the $\min \{g_1,g_2,\dots, g_l \}>1,$ number 1 cannot be obtained in the form
\begin{eqnarray*}
\sum_{i=1}^l \lambda_i g_i =1,\;\;\;\forall i:\,\lambda_i\in \mathbb{Z}_+.
\end{eqnarray*}

In general, a cone is called pointed if the origin is its vertex. It is equivalent to the property that the origin is not an interior point
of segment such that it is part of the cone.
In what follows, the unique minimal Hilbert basis of a pointed cone
related to the FCC grid is uncovered. It is done by a method which works at
other cones as well. An important notion used by the method is the {\em support of a non-negative vector}. It is the set of the indices of the positive
components. The notation is as follows: if $x^T=(x_1,x_2,\dots,x_n)$, then $supp(x)=\{i\mid x_i>0\}$. The size of $supp(x)$ is the number of its elements
denoted by $\mid supp(x)\mid \,=\,\mid \{i\mid x_i>0\}\mid$. There are cases such that if a vector has a minimal support, then the vector belongs to the
minimal Hilbert basis.

\begin{theorem}\label{THB3}
Let $C$ be a pointed rational polyhedral cone of non-negative vectors and $K=C\cap \mathbb{Z}^n$. Assume that $H$ is the unique minimal Hilbert basis of $C$.
Let $L$ be a subset of $H$.
Assume that the vector $x^*$
satisfies the conditions as follows:
\begin{eqnarray}\label{No1}
x^*\in argmin\left\{ \mid supp(x)\mid : x\in K\setminus L\right\}
\end{eqnarray}
and
\begin{eqnarray}\label{No2}
\forall x\in L:\; supp(x)\not\subset supp(x^*)
\end{eqnarray}
and
\begin{eqnarray}\label{No3}
0\not=x^*\leq e.
\end{eqnarray}
Then $x^*\in H$.
\end{theorem}
\proof{
If $x^*\not\in H$, then there are two elements $u$ and $v$ of $H$ such that $u,\,v\not= 0$ and $u+v\leq x^*$. It follows from (\ref{No3}) that
$\forall i:\,u_i v_i=0$. If $v\not\in L$, then $v$ has a smaller support size than $x^*$, that is $x^*$ cannot be an optimal solution in (\ref{No1}).
Thus, $u,\,v\in L$ violates (\ref{No2}).\hfill Q.E.D.\\
}

The next theorem is similar to the last one. A significant difference is that the optimization problem of (\ref{No1}) is substituted by lexicographic
optimization of two objective functions.

\begin{theorem}\label{THB3.5}
Let $C$ be a pointed  polyhedral cone of non-negative vectors and $K=C\cap \mathbb{Z}^n$. Assume that $H$ is the unique minimal Hilbert basis of $C$.
Let $L$ be a subset of $H$. Assume that the vector $x^*$
satisfies the conditions as follows. It is the lexicographic optimal solution of the problem
\begin{eqnarray}\label{No4}
l-min\left\{ \bigg(\begin{array}{l}\mid supp(x)\mid\\ \sum_{i=1}^{n}x_i \end{array} \bigg): x\in K\setminus L\right\}
\end{eqnarray}
and
\begin{eqnarray}\label{No5}
\forall x\in L:\; supp(x)\not\subset supp(x^*)
\end{eqnarray}
and
\begin{eqnarray}\label{No6}
0\not=x^*\leq 2e.
\end{eqnarray}
Then $x^*\in H$.
\end{theorem}
\proof{
First, assume that there is a vector $y\in\{0,1\}^n$ such that $x^*=2y$. Notice that $\mid supp(x^*)\mid=\mid supp(y)\mid$. If $y\not\in L$, then $y$
is a lexicographically better solution than $x^*$. Thus, $x^*$ is not optimal in (\ref{No4}). Hence, $y\in L$. However, (\ref{No5}) is violated in this case. Thus, $x^*$ has no form that  $x^*=2y$. If $x^*\not\in H$, then it must be the sum of some elements of $H$.
If the sum cannot consists of a single element which is $x^*$ itself, then $x^*\in H$.
Assume that for all $u$ and $v$ of the sum: $supp(u)\cap supp(v) =\emptyset$.
Then either all elements of the sum are the elements of $L$ as well which means that (\ref{No5}) is violated, or at least one of the elements of
the sum is not in $L$ implying that $x^*$ is not optimal in (\ref{No4}).
The only case remained that there are $u$ and $v$ in the sum, with
$supp(u)\cap supp(v)\not =\emptyset$, $supp(v)\not\subset supp(u)$ and $supp(u)\not\subset supp(v)$. Hence again, if all elements of the sum are the elements of $L$,
then (\ref{No5}) is violated. Otherwise, $x^*$ is not optimal in (\ref{No4}).
\hfill Q.E.D.\\
}

Let $A'$ be an $m \times n$ rational matrix. The cone is defined in the form
\begin{eqnarray}\label{HB05}
C=\left\{x\in \mathbb{R}_+^n\mid A'x=0\right\},
\end{eqnarray}
where $0$ is the zero vector. It is obvious that $C$ satisfies the algebraic properties.

\begin{theorem}\label{TH4}
If $C$ satisfies (\ref{HB05}) and $v\leq e$ is an element of $H$ the unique minimal Hilbert basis of $K$, then there is no $u\in H$ for which $supp(v)\subseteq supp(u)$.
\end{theorem}
\proof{
Let $u$ be an element of $K$ for which $supp(v)\subseteq supp(u)$.
Let $k$ be the minimal non-zero component of $u$ ($k$ is positive).
Then the vector $u-kv$ is also an element of $K$, because $A'(u-kv)=0$ and the elements of $u-kv$ are non-negative.
Thus, $u-kv\in K$. Hence for some non-negative integer coefficients
\begin{eqnarray}\label{HB06}
u-kv=\sum_{h\in H}\lambda_h h,
\end{eqnarray}
where $\forall h: \lambda_h\in Z_{+}$. Now $\lambda_v=0$ and if $u\in H$, then $\lambda_u=0$.
Thus
\begin{eqnarray}\label{HB07}
u=\sum_{h\in H\setminus\{v\}}\lambda_h h + kv.
\end{eqnarray}
In (\ref{HB07}), the coefficient of $v$ is positive and $v\not =u$.
It means that $u$ cannot be an element of the minimal Hilbert basis of $K$.
\hfill Q.E.D.\\}

Note that in case of the FCC grid matrix $A$ will be used in the role of $A'$ in the previous Theorem.
\section{The minimal Hilbert Basis of FCC}\label{HB}

Whenever we consider the Hilbert Basis for the FCC grid, we consider 18 dimensional vectors $w=(w_1,w_2,\ldots,w_{18})$,  which satisfies $A w=0$. Here, $w_j \in \mathbb Z_+$, $j=1,2,\ldots,18$ where such a vector describes how many times step $i$ is used. The space 
of these vectors will be referred as the ``step-space'', and elements of the
a Hilbert basis are 
 vectors satisfying the above condition of the step-space.

The elements of the minimal Hilbert basis of the cone of the FCC grid are provided in Table 1. The example $a_1+a_4=0$ in the first row
of the table means the vector of the step-space as follows:
\begin{eqnarray*}
(1,0,0,1,0,0,0,0,0,0,0,0,0,0,0,0,0,0).
\end{eqnarray*}
Similarly, the example of type 3 is described by the vector
\begin{eqnarray*}
(1,0,0,0,0,0,0,0,1,1,0,0,0,0,0,0,0,0).
\end{eqnarray*}

\begin{table}[H]
\caption{Types of the elements of the minimal Hilbert basis of the FCC grid.}\label{tab1}
\begin{center}
\begin{tabular}{|c|c|c|c|c|}
\hline
Type & \# of $a$'s & \# of $d$'s & Example & \# of cases\\ \hline
1 & 2 & 0 & $a_1+a_4=0$ & 3\\
2 & 0 & 2 & $d_1+d_{4}=0$ & 6\\
3 & 1 & 2 & $a_1+d_{3}+d_{4}=0$ & 12\\
4 & 0 & 3 & $d_1+d_{7}+d_{12}=0$ & 8\\
5 & 1 & 3 & $a_1+d_{4}+d_{8}+d_{9}=0$ & 24\\
6 & 0 & 4 & $d_1+d_2+d_7+d_{8}=0$ & 6\\
7 & 2 & 3 & $a_4+a_5+d_{1}+d_{5}+d_{10}=0$ & 24\\
8 & 2 & 2 & $a_4+a_5+2d_1=0$ & 12\\
9 & 3 & 3 & $a_1+a_2+a_3+d_{4}+d_{8}+d_{12}=0$ & 8\\
10 & 1 & 4 & $a_1+2d_3+d_{11}+d_{12}=0$ & 24\\
11 & 2 & 4 & $a_1+a_2+2d_7+2d_{12}=0$ & 24\\
\hline
\end{tabular}
\end{center}

\end{table}

The remainder of this section consists of three parts. In the first subsection, important properties of the vectors of the step-space are analyzed.
In the second subsection, it is shown that Table 1 contains all elements of the minimal Hilbert basis and contains no superfluous elements.
Finally, in the third subsection, geometric interpretation is presented. Moreover,
it is shown that of each type of elements there are exactly as many as are indicated in the last column of the table.

\subsection{Some properties of the vectors of the step-space}

\begin{theorem}\label{TH5}
If $supp(v)\subseteq supp(u)$, where $v$ is an element of Types 1-7 and 9 in Table \ref{tab1}, then $u$ is not element of $H$, the minimal Hilbert basis of $K$.
\end{theorem}
\proof{For all elements of Types 1-7 and 9 $v\leq e$ holds, thus this statement is the corollary of Theorem \ref{TH4}.\hfill Q.E.D.\\}

\begin{theorem}\label{TH6}
If there is a component of $u$ belonging to $a_i$, which value is at least 2, i.e. $\max\{u_1, u_2, \dots, u_6\}\geq 2$, then $u$ is not element of $H$, the minimal Hilbert basis of $K$.
\end{theorem}
\proof{For example let $u_1\geq2$. Now $u\in K$, i.e. $Au=0$, where $A$ is the matrix of (\ref{e001}). It means that
\begin{eqnarray}\label{HB08}
a_1u_1+a_2u_2+\dots, +a_6u_6+d_1u_7+d_2u_8+\dots+d_{12}u_{18}=0.
\end{eqnarray}
If the component $u_1$ belonging to $a_1=(2,0,0)^T$ is at least 2, then the first coordinate of $a_1u_1$ is at least 4.
The right hand side of (\ref{HB08}) is 0 also in the first coordinate.
For reaching this, $u$ needs to have positive components for the vectors, where the sign of the first coordinate is negative: $a_4, d_3, d_4, d_7, d_8$.

If $u_1\geq2$, and $u_4$ is positive, then $a_1+a_4=0$ is part of $u$, thus $u$ is not element of $H$.

Let $u_1\geq2$, $u_4=0$.
If $u_9,u_{10}\geq1$, then $a_1+d_3+d_4=0$ is part of $u$. Thus, $u$ is not an element of $H$.
If at least one of $u_9,u_{10}$ is equal to 0, then at least one of the components belonging to $d_3, d_4, d_7, d_8$ (where the first coordinate is $-1$) is at least 2.
For example let $u_9\geq 2$. Then $u_{10}=0$.
Moreover the value of the second coordinate of $u_9d_3\geq2$.
For reaching 0, $u$ needs to have positive components for the vectors, where the sign of the second coordinate is negative: $a_5, d_2, d_4, d_{11}, d_{12}$.

Now $u_1,u_9\geq2$, $u_4,u_{10}=0$.
If $u_5\geq1$, then $a_1+a_5+2d_3=0$ is part of $u$. Thus, $u$ is not an element of $H$.
If $u_8\geq1$, then $d_2+d_3=0$ is part of $u$. Thus, $u$ is not an element of $H$.
Only $d_{11}, d_{12}$ left from the group of vectors where the sign of the second coordinate is negative.
If both of the components belonging to them ($u_{17}$ and $u_{18}$) are at least one, then $a_1+2d_3+d_{11}+d_{12}=0$ is part of $u$. Thus, $u$ is not an element of $H$.
Then for example let $u_{17}\geq 2$ and $u_{18}=0$.
In this case the value of the third coordinate of $u_{17}d_{11}\geq2$.
For reaching 0, $u$ needs to have positive components for the vectors, where the sign of the third coordinate is negative: $a_6, d_6, d_8, d_{10}, d_{12}$.

Now $u_1,u_9,u_{17}\geq2$, $u_4,u_5,u_8,u_{10},u_{18}=0$.
If $u_{12}\geq1$, then $d_3+d_6+d_{11}=0$ is part of $u$. Thus, $u$ is not an element of $H$.
If $u_{14}\geq1$, then $a_1+d_3+d_8+d_{11}=0$ is part of $u$. Thus, $u$ is not an element of $H$.
If $u_{16}\geq1$, then $d_{10}+d_{11}=0$ is part of $u$. Thus, $u$ is not an element of $H$.
No more vectors, where the sign of the third coordinate is negative. Thus, $Au=0$ is not true for $u$.

\hfill Q.E.D.\\}

\begin{theorem}\label{TH7}
If there is a positive component of $u$ belonging to $a_i$ and there is another positive component belonging to $d_j$, where $a_i$ and $d_j$ have the same sign in a given coordinate, then $u$ is not an element of the minimal Hilbert basis of $K$.
\end{theorem}
\proof{
Let $H$ be again the minimal Hilbert basis.
For example let the component $u_1\geq1$ (belonging to vector $a_1$) and let the component $u_7\geq1$ (belonging to $d_1$).
Then the first coordinate of $a_1u_1+d_1u_7$ is at least 3.
The right hand side of (\ref{HB08}) is 0 in the first coordinate.
For reaching this, $u$ needs to have positive components for the vectors, where the sign of the first coordinate is negative: $a_4, d_3, d_4, d_7, d_8$.

If $u_4$ is positive, then $a_1+a_4=0$ is part of $u$. Thus, $u$ is not an element of $H$.
If $u_{10}$ is positive, then $d_1+d_4=0$ is part of $u$. Thus, $u$ is not an element of $H$.
Only $d_3, d_7, d_8$ left from the group of vectors where the sign of the first coordinate is negative.

\underline{Case 1}
Let the component $u_9\geq1$ (belonging to $d_3$).
Now $u_1,u_7,u_9\geq1$ and $u_4,u_{10}=0$.
The value of the second coordinate of the vector $d_1u_7+d_3u_9$ is at least 2.
For reaching 0, $u$ needs to have positive components for the vectors, where the sign of the second coordinate is negative: $a_5, d_2, d_4, d_{11}, d_{12}$ (but $u_{10}=0$ now).

If $u_5\geq1$, then $a_5+d_1+d_3=0$ is part of $u$. Thus, $u$ is not an element of $H$.
If $u_8\geq1$, then $d_2+d_3=0$ is part of $u$. Thus, $u$ is not an element of $H$.
Only $d_{11}, d_{12}$ left from the group of vectors where the sign of the second coordinate is negative.
If both of the components belonging to them ($u_{17}$ and $u_{18}$) are at least one, then $d_1+d_3+d_{11}+d_{12}=0$ is part of $u$. Thus, $u$ is not an element of $H$.
Assume for example that $u_{17}\geq 2$ and $u_{18}=0$.
In this case the value of the third coordinate of $u_{17}d_{11}\geq2$.
For reaching 0, $u$ needs to have positive components for the vectors, where the sign of the third coordinate is negative: $a_6, d_6, d_8, d_{10}, d_{12}$ (but $u_{18}=0$ now).

According to the current assumption $u_1,u_7,u_9\geq1$ and $u_{17}\geq2$.
If $u_6\geq1$, then $a_6+d_1+d_3+2d_{11}=0$ is part of $u$. Thus, $u$ is not an element of $H$.
If $u_{12}\geq1$, then $d_3+d_6+d_{11}=0$ is part of $u$. Thus, $u$ is not an element of $H$.
If $u_{14}\geq1$, then $d_1+d_8+d_{11}=0$ is part of $u$. Thus, $u$ is not an element of $H$.
If $u_{16}\geq1$, then $d_{10}+d_{11}=0$ is part of $u$. Thus, $u$ is not an element of $H$.
There are no more vectors, where the sign of the third coordinate is negative. Thus, $Au=0$ is not true for $u$.

\underline{Case 2}
Assume that $u_9=0$. Assume also that $u_{13}\geq2$ (belonging to $d_7$).
Now $u_1,u_7\geq1$, $u_{13}\geq2$ and $u_4,u_9,u_{10}=0$.
The value of the third coordinate of $d_7u_{13}$ is at least 2.
For reaching 0, $u$ needs to have positive components for the vectors, where the third coordinate is negative: $a_6, d_6, d_8, d_{10}, d_{12}$.

If $u_6\geq1$, then $a_1+a_6+2d_7=0$ is part of $u$. Thus, $u$ is not an element of $H$.
If $u_{12}\geq1$, then $d_6+d_7=0$ is part of $u$. Thus, $u$ is not an element of $H$.
If $u_{14}\geq1$, then $a_1+d_7+d_8=0$ is part of $u$. Thus, $u$ is not an element of $H$.
If $u_{18}\geq1$, then $d_1+d_7+d_{12}=0$ is part of $u$. Thus, $u$ is not an element of $H$.
Only $d_{10}$ left from the group of vectors where the sign of the third coordinate is negative.
It means that $u_{16}\geq2$ holds.

The value of the second coordinate of $d_1u_7+d_{10}u_{16}$ is at least 3.
For reaching 0, $u$ needs to have positive components for the vectors, where the sign of the second coordinate is negative: $a_5, d_2, d_4, d_{11}, d_{12}$ (but $u_{10}=u_{18}=0$ now).

Now $u_1,u_7\geq1$, $u_{13},u_{16}\geq2$.
If $u_5\geq1$, then $a_1+a_5+2d_7+2d_{10}=0$ is part of $u$. Thus, $u$ is not an element of $H$.
If $u_8\geq1$, then $d_2+d_7+d_{10}=0$ is part of $u$. Thus, $u$ is not an element of $H$.
If $u_{17}\geq1$, then $d_{10}+d_{11}=0$ is part of $u$. Thus, $u$ is not an element of $H$.
No more vectors, where the sign of the second coordinate is negative. Thus, $Au=0$ is not true for $u$.
\hfill Q.E.D.\\}

\begin{theorem}\label{TH8}
If for a given $u$ the sum of the components belonging to $d_j$'s, which have the same sign for a coordinate, is at least 3, then $u$ is not an element of $H$, the minimal Hilbert basis of $K$.
\end{theorem}

\proof{
Denote the minimal Hilbert basis by $H$ again. Assume that the coordinate in question is the first one
and let the sign be positive.
Then the vectors, for which the first coordinate is positive: $d_1, d_2, d_5, d_6$.

\underline{Case 1}
Let the component $u_7\geq3$ (belonging to vector $d_1$).
Then the first coordinate of $d_1u_7$ is at least 3.
The right hand side of (\ref{HB08}) is 0 in the first coordinate.
For reaching this, $u$ needs to have positive components for the vectors, where the sign of the first coordinate is negative: $a_4, d_3, d_4, d_7, d_8$.

If $u_4$ is positive, then $u$ is not an element of $H$ because of the previous two theorems.
If $u_{10}$ is positive, then $d_1+d_4=0$ is part of $u$. Thus, $u$ is not an element of $H$.
Only $d_3, d_7, d_8$ left from the group of vectors where the sign of the first coordinate is negative.

\underline{Case 1a}
Let the component $u_{13}\geq1$ (belonging to $d_7$).
The value of the second coordinate of $d_1u_7$ is at least 3.
For reaching 0, $u$ needs to have positive components for the vectors, where the sign of the second coordinate is negative: $a_5, d_2, d_4, d_{11}, d_{12}$ (remember that $u_{10}=0$ now).

If $u_5\geq1$, then $u$ is not an element of $H$, because of the previous two theorems.
If $u_{18}\geq1$, then $d_1+d_7+d_{12}=0$ is part of $u$. Thus, $u$ is not an element of $H$.
Only $d_2, d_{11}$ left from the group of vectors where the sign of the second coordinate is negative.

\underline{Case 1a1}
Let $u_8\geq1$ (belonging to $d_2$).
The value of the third coordinate of $u_{13}d_7\geq1$.
For reaching 0, $u$ needs to have positive components for the vectors, where the sign of the third coordinate is negative: $a_6, d_6, d_8, d_{10}, d_{12}$ (remember that $u_{18}=0$ now).

Now $u_7\geq3$ and $u_8,u_{13}\geq1$.
If $u_{12}\geq1$, then $d_6+d_7=0$ is part of $u$. Thus, $u$ is not an element of $H$.
If $u_{14}\geq1$, then $d_1+d_2+d_7+d_8=0$ is part of $u$. Thus, $u$ is not an element of $H$.
If $u_{16}\geq1$, then $d_2+d_7+d_{10}=0$ is part of $u$. Thus, $u$ is not an element of $H$.
It means that for reaching 0 in the third coordinate: $u_6\geq1$ (belonging to $a_6$).
But $u_6$ cannot be at least 2 because of Theorem \ref{TH6}. Thus, $u_6=1$.
Now the sum of third coordinates, for the vectors where the third coordinate is negative, is equal to $-2$.
For reaching 0, $u$ needs to have positive components with the sum equal to 2, for the vectors, where the sign of the third coordinate
is positive: $a_3, d_5, d_7, d_9, d_{11}$.
Now the sum of the first two coordinates of $d_1u_7+d_2u_8$ is at least 6.
The set $a_3, d_5, d_7, d_9, d_{11}$ can decrease this 6 only by 2, and the remaining vectors, where the third coordinate is equal to 0 ($d_2, d_3$) cannot decrease this sum.
It means the reaching 0 in the first two coordinates is impossible.

\underline{Case 1a2}
Let $u_{17}\geq1$ (belonging to $d_{11}$).
Now $u_7\geq3$ and $u_{13},u_{17}\geq1$.
The value of the third coordinate of $u_{13}d_7+u_{17}d_{11}\geq2$.
For reaching 0, $u$ needs to have positive components for the vectors, where the sign of the third coordinate is negative: $a_6, d_6, d_8, d_{10}, d_{12}$ (but $u_{18}=0$ now).
If $u_6\geq1$, then $a_6+d_1+d_7+d_{11}=0$ is part of $u$. Thus, $u$ is not an element of $H$.
If $u_{12}\geq1$, then $d_6+d_7=0$ is part of $u$. Thus, $u$ is not an element of $H$.
If $u_{14}\geq1$, then $d_1+d_8+d_{11}=0$ is part of $u$. Thus, $u$ is not an element of $H$.
If $u_{16}\geq1$, then $d_{10}+d_{11}=0$ is part of $u$. Thus, $u$ is not an element of $H$.
Thus $Au=0$ is not true for $u$.

\underline{Case 1b}
The case when the component $u_{14}\geq1$ (belonging to $d_8$) is the same as Case 1a, only the sign of the third coordinate is the opposite.

\underline{Case 1c}
If $u_{13}=u_{14}=0$, then for reaching 0, $u_9\geq1$ (belonging to $d_3$).
Now $u_7\geq3$ and $u_9\geq1$. Thus, the second coordinate of $d_1u_7+d_3u_9$ is at least 4.
For reaching 0, $u$ needs to have positive components for the vectors, where the sign of the second coordinate is negative: $a_5, d_2, d_4, d_{11}, d_{12}$ (but $u_{10}=0$ now).
If $u_5\geq1$, then $u$ is not an element of $H$, because of the previous two theorems.
If $u_8\geq1$, then $d_2+d_3=0$ is part of $u$. Thus, $u$ is not an element of $H$.
Only $d_{11},d_{12}$ left from the group of vectors where the sign of the second coordinate is negative.

For example if $u_{17}\geq1$, then the third coordinate of $u_{17}d_{11}\geq1$.
For reaching 0, $u$ needs to have positive components for the vectors, where the sign of the third coordinate is negative: $a_6, d_6, d_8, d_{10}, d_{12}$ (but $u_{14}=0$ now).
If $u_{12}\geq1$, then $d_3+d_6+d_{11}=0$ is part of $u$. Thus, $u$ is not an element of $H$.
If $u_{16}\geq1$, then $d_{10}+d_{11}=0$ is part of $u$. Thus, $u$ is not an element of $H$.
If $u_{18}\geq1$, then $d_1+d_3+d_{11}+d_{12}=0$ is part of $u$. Thus, $u$ is not an element of $H$.
It means that for reaching 0 in the third coordinate: $u_6\geq1$ (belonging to $a_6$).
Notice that this situation is the same as in Case 1a1.

\underline{Case 2}
Let $u_7=2$ (belonging to $d_1$) and let another component be at least 1 belonging to $d_2,d_5,d_6$.
Then the first coordinate of $d_1u_7$ is at least 3, again.
For reaching 0, $u$ needs to have positive components for the vectors, where the sign of the first coordinate is negative: $a_4, d_3, d_4, d_7, d_8$.

\underline{Case 2a}
Let $u_7=2$ and $u_8\geq1$ (belonging to $d_2$).
If $u_4$ is positive, then $u$ is not an element of $H$ because of the previous two theorems.
If $u_9$ is positive, then $d_2+d_3=0$ is part of $u$. Thus $u$ is not an element of $H$.
If $u_{10}$ is positive, then $d_1+d_4=0$ is part of $u$. Thus, $u$ is not an element of $H$.
If $u_{13}$ and $u_{14}$ are positive, then $d_1+d_2+d_7+d_8=0$ is part of $u$. Thus, $u$ is not an element of $H$.
But if only one of $u_{13}$ and $u_{14}$ is positive, then its value is at least 3.
This is Case 1.

\underline{Case 2b}
Let $u_7=2$ and $u_{11}\geq1$ (belonging to $d_5$).
(The case when $u_{12}\geq1$ is the same, only this sign of the third coordinate is different.)
If $u_4$ is positive, then $u$ is not an element of $H$ because of the previous two theorems.
If $u_{10}$ is positive, then $d_1+d_4=0$ is part of $u$. Thus, $u$ is not an element of $H$.
If $u_{14}$ is positive, then $d_5+d_8=0$ is part of $u$. Thus, $u$ is not an element of $H$.
If only one of $u_9$ and $u_{10}$ is positive, then its value is at least 3, and this is Case 1.

Let $u_9$ and $u_{10}$ positive.
Now $u_7=2$ and $u_9,u_{11},u_{12}\geq1$. Thus, the second coordinate of $d_1u_7+d_3u_9$ is at least 3.
For reaching 0, $u$ needs to have positive components for the vectors, where the sign of the second coordinate is negative: $a_5, d_2, d_4, d_{11}, d_{12}$ (but $u_{10}=0$ now).
If $u_5\geq1$, then $u$ is not an element of $H$, because of the previous two theorems.
If $u_8\geq1$, then $d_2+d_3=0$ is part of $u$. Thus, $u$ is not an element of $H$.
If $u_{18}$ is positive, then $d_1+d_7+d_{12}=0$ is part of $u$. Thus, $u$ is not an element of $H$.
Only $d_{11}$ left.

Let $u_{17}\geq1$, then the third coordinate of $u_{11}d_5+u_{13}d_7+u_{17}d_{11}\geq3$.
For reaching 0, $u$ needs to have positive components for the vectors, where the sign of the third coordinate is negative: $a_6, d_6, d_8, d_{10}, d_{12}$ (but $u_{18}=0$ now).
If $u_6\geq1$, then $u$ is not an element of $H$, because of the previous two theorems.
If $u_{12}\geq1$, then $d_6+d_7=0$ is part of $u$. Thus, $u$ is not an element of $H$.
If $u_{14}\geq1$, then $d_5+d_8=0$ is part of $u$. Thus, $u$ is not an element of $H$.
If $u_{16}\geq1$, then $d_{10}+d_{11}=0$ is part of $u$. Thus, $u$ is not an element of $H$.
Thus, $Au=0$ is not true for $u$.

\underline{Case 3}
Let $u_7=1$ (belonging to $d_1$) and let another two components be 1 belonging to $d_2,d_5,d_6$, for example $u_8$ and $u_{11}$.
Then the first coordinate is at least 3, again.
For reaching 0, $u$ needs to have positive components for the vectors, where the sign of the first coordinate is negative: $a_4, d_3, d_4, d_7, d_8$.
If $u_4$ is positive, then $u$ is not an element of $H$ because of the previous two theorems.
If $u_9$ is positive, then $d_2+d_3=0$ is part of $u$. Thus, $u$ is not an element of $H$.
If $u_{10}$ is positive, then $d_1+d_4=0$ is part of $u$. Thus, $u$ is not an element of $H$.
If $u_{14}$ is positive, then $d_5+d_8=0$ is part of $u$. Thus, $u$ is not an element of $H$.
Only $d_7$ left, but in this case its component value is at least 3.
This is Case 1.
\hfill Q.E.D.\\}

\subsection{The completeness and exactness of Table 1}

The completeness means that Table 1 contains all elements of the minimal Hilbert basis. The exactness means that
the minimal Hilbert basis does not contain any other element.

The logic of the argument is that any solution such that for some reason, cannot be an element of the minimal Hilbert basis should be excluded.
For this, linear conditions are defined. The elements of the step-space are denoted by $(x_1,\dots,x_{18})^T\in\mathbb{Z}^{18}$. $A$ is still
the matrix of (\ref{e001}). The first constraints are the equations of
\begin{eqnarray}\label{e501}
Ax=0.
\end{eqnarray}
The zero vector is not an element of the Hilbert basis. That is why the inequality
\begin{eqnarray}\label{e502}
\sum_{j=1}^{18} x_j \geq 1
\end{eqnarray}
is claimed.

What Theorems \ref{TH6}-\ref{TH8} say on the first row of matrix $A$ is that if the terms with positive coefficients
are summed up, then the sum may not exceed 2, i.e.
\begin{eqnarray}\label{e503}
2x_1+x_7+x_8+x_{11}+x_{12} \leq 2.
\end{eqnarray}
A similar inequality holds for the negative terms
\begin{eqnarray}\label{e504}
2x_4+x_9+x_{10}+x_{13}+x_{14} \leq 2.
\end{eqnarray}
The theorems are valid to the second and third rows of matrix $A$ as well, giving four more inequalities, as follows:
\begin{eqnarray}\label{e505}
2x_2+x_7+x_9+x_{15}+x_{16} \leq 2\\
2x_5+x_8+x_{10}+x_{17}+x_{18} \leq 2\\
2x_3+x_{11}+x_{13}+x_{15}+x_{17} \leq 2\\
2x_6+x_{12}+x_{14}+x_{16}+x_{18} \leq 2.
\end{eqnarray}

A large part of the investigation is based on the support sets of the solutions listed in Table 1.
Theorem \ref{TH5} implies that all solutions having a support set
at least as large as an element of Types 1-7 or 9 can be excluded. Further variables are needed to describe the support sets.
They are binary indicator variables showing if an index is included in support set. These variables are denoted by $y_7,y_8,\dots,y_{18}$. Theorem
\ref{TH6} implies that the values of the variables $x_1,\dots,x_6$ can be only either 0, or 1, i.e. these variables are themselves binary.
Thus, they can be used as their own indicator variables. The components of the vectors of the step-space and the relevant indicator variables
must be related in a way that a component is positive if and only if the value of the indicator variable is 1. It implies that
\begin{eqnarray}\label{e506}
y_j \leq x_j,\;\,j=7,\dots,18.
\end{eqnarray}
Hence, if $x_j=0$, then $y_j=0$, too. The opposite direction can be claimed by the so-called big-$M$ method of optimization. Let $M$ be a very large
positive number. Then the inequality
\begin{eqnarray}\label{e507}
x_j \leq My_j,\;\,j=7,\dots,18
\end{eqnarray}
claims that if $x_j>0$, then $y_j=1$ as $y_j$ is restricted either 0, or 1.

The constraint that the support set of a solution is not at least as large as
the support set of a solution listed in Table 1, is equivalent to the linear inequality that the sum of the indicator variables of the potential
new solution on the indices of the listed solution is less than the size of the support set of the listed solution. It means in the case of the
example of Type 1 in the table that
\begin{eqnarray*}
x_1+x_4 \leq 1.
\end{eqnarray*}
The most complicated case is that of Type 9. The relevant inequality in the case of the example in the table is as follows
\begin{eqnarray*}
x_1 +x_2 +x_3 +y_{10} +y_{14}+y_{18} \leq 5.
\end{eqnarray*}
The total number of these kinds of inequalities is 91 according to the last column of Table 1. In what follows, these inequalities are referred as
Type-179 inequalities.

Types 8, 10, and 11 need more careful investigation. The step-space vectors have a component with value 2.
In the case, if this component is forced to be less, other solutions may exist. Each case is investigated in two steps. First, it is shown, that the
solution is unique on the same support set under the above mentioned constraints. The second step is proving that if a larger support set is allowed,
but at least one of the components having value 2 is restricted to be less than 2, then there is no solution.

{\em The uniqueness of the solution of Type 8}

Only three variables are allowed to be different from 0. The example of Type 8 in Table 1 is $a_4+a_5+2d_1$
giving a small equation system as follows:
\begin{eqnarray*}
-2x_4 + x_7 = 0\\ -2x_5 + x_7 = 0.
\end{eqnarray*}
There is no third equation as in the last row all coefficients are 0. It is easy to see that in every solution of the two equations
$2x_4=2x_5=x_7$. Hence, the only non-negative and non-zero integer solution without a variable having a value greater than 2 is $x_4=x_5=1, \,x_7=2$.
As all solutions of Type 8 have the same structure, the same method works for all of them.

{\em The uniqueness of the solution of Type 10}

The example for Type 10 is $a_1+2d_3+d_{11}+d_{12}$. The relevant equation system is as follows:
\begin{eqnarray*}
2x_1-x_9 = 0\\ x_9-x_{17}-x_{18} = 0\\ x_{17}-x_{18} = 0.
\end{eqnarray*}
Hence, the relations $x_9=2x_1,\,x_{17}=x_{18},\,x_9=2x_{17}$ are obtained. They imply that $x_{17}=x_{18}=x_1$. Thus, the form of all solutions
restricted to the indices 1, 9, 17, and 18 is $(w,2w,w,w)$. Notice that anything else is 0. Similarly to the previous case,
the only non-negative and non-zero integer solution without a variable having a value greater than 2 is $x_1=1,\,x_9=2,\,x_{17}=1,$ and $x_{18}=1$.

{\em The uniqueness of the solution of Type 11}

The example for Type 11 is $a_1+a_2+2d_7+2d_{12}$. The relevant equation system is as follows:
\begin{eqnarray*}
2x_1-x_{13} = 0\\ 2x_2-x_{18} = 0\\ x_{13}-x_{18} = 0.
\end{eqnarray*}
Hence, the relations $x_{13}=2x_1,\,x_{18}=2x_{2}\,x_{13}=x_{18},$ are obtained. They imply that $x_1=x_2$. Thus, the form of all solutions
restricted to the indices 1, 2, 13, and 18 is $(w,w,2w,2w)$. Notice that anything else is 0. Similarly to the previous case,
the only non-negative and non-zero integer solution without a variable having a value greater than 2 is $x_1=1,\,x_2=1,\,x_{13}=2,$ and $x_{18}=2$.

Now, the second step is to be done. It proves that if a larger support set is allowed,
but at least one of the components having value 2 is restricted to be less than 2, then there is no solution. The technique of showing non-existence
is as follows. A constraint set is defined. It consists of constraints (\ref{e501})-(\ref{e507}) and the Type-179 inequalities. This constraint set is
the same at the investigations of Types 8, 10, and 11. Further constraints are added according to the particular type. The types are discussed with the
same member of them.

{\em No solution exists with larger support set in case of Type 8}

The example of Type 8 in Table 1 is $a_4+a_5+2d_1$. Thus, variables $x_4,\,x_5$, and $x_7$ must be positive. Theorem \ref{TH6} implies that the only
possible value of $x_4,$ and $x_5$ is 1. Variable $x_7$ must be positive and less than 2. Thus, its only possible value is 1. Hence, the further constraints
are
\begin{eqnarray}\label{e508}
x_4=1,\,x_5=1,\,x_7=1.
\end{eqnarray}

{\em No solution exists with larger support set in case of Type 10}

The example of Type 10 in Table 1 is $a_1+2d_3+d_{11}+d_{12}$. The same logic applies for variables $x_1$, and $x_9$. Variables $x_{17}$, and $x_{18}$ can be
greater than 1. Thus, the further constraints are
\begin{eqnarray}\label{e509}
x_1=1,\,x_9=1,\,x_{17}\geq 1,\,x_{18}\geq 1.
\end{eqnarray}

It is worth to finish these two types before discussing Type 11. There are two constraint sets as follows:
\begin{center}
(\ref{e501})-(\ref{e507}), Type-179 inequalities, and (\ref{e508})
\end{center}
and
\begin{center}
(\ref{e501})-(\ref{e507}), Type-179 inequalities, and (\ref{e509}).
\end{center}
Each is the constraint set of an integer programming problem in the optimization theory. Thus, they can be investigated by an optimization program
if any objective function is added. LINGO 18.0 was applied in that case. The result for both types was that no feasible solution exists. The conclusion
what one can obtain is that the support sets can be excluded. This conclusion is still not equivalent to the exactness of the cases contained in Table 1.

{\em No solution exists with larger support set in case of Type 11}

The example of Type 11 in Table 1 is $a_1+a_2+2d_7+2d_{12}$. If a similar logic is applied, then the constraints $x_1=1$, $x_2=1$, and $x_{13}+x_{18}\leq 3$
should be added. The explanation of the latter is that the variables may not be greater than 2 and if both of them are 2, then the same solution is
obtained again. When only these constraints are introduced, then the optimizer software gives the feasible solution $a_1+a_2+2d_8+2d_{11}$ which is
the pair of the first one as $d_7+d_{12}=d_8+d_{11}$. As a matter of fact, $a_1+a_2+2d_8+2d_{11}$ is another element of the Hilbert basis belonging to
Type 11. Thus, it also must be excluded by the constraint $x_{14}+x_{17}\leq 3$. After adding this fourth constraint, no feasible solution remained.
The conclusion is that the two support sets can be excluded by a single optimization.

The constraint sets excluding the supports sets of Types 8, 10, and 11 are denoted by Type-8, Type-10, and Type-11 constraints.

The final evidence for the completeness of the Hilbert basis of Table 1 is obtained by the optimizer. Another constraint set is formed. It includes
\begin{center}
{\em constraints (\ref{e501})-(\ref{e507}) and the Type-179 inequalities, Type-8, Type-10, and Type-11 constraints and the integrality of the variables.}
\end{center}
It is again the constraint set of an integer programming problem. The optimizer shows that it has no feasible solution. The meaning of the results is
that no feasible solution of the equation system exists such that it satisfies the restrictions of the Theorems \ref{TH5}-\ref{TH8} and its support set
is different from the support sets of the known elements of the Hilbert basis.

\subsection{The Geometric Meaning of the Elements of the FCC Hilbert Basis}
An element  $w$ of the Hilbert basis is a minimal vector of the step-space, i.e., solution of $Aw=0$. In other words we cannot find another eighteen dimensional nonzero vector $u$ in step-space such that $u\leq w$ and $Au=0$. The vector $w$ can be defined as sum of the vectors giving zero-sum in grid-space.
This gives a geometric interpretation: an element of the Hilbert basis can be interpreted as a directed cycle in the grid that cannot be decomposed to smaller cycles (in terms of less number of steps).

Figures \ref{type1}, \ref{type4},  \ref{type7}, and \ref{type10} represent different types of the Hilbert basis elements. Before giving the details of these figures, we fix the coordinate system: in each figure, the point at the bottom-left corner represents the origin, i.e. $(0, 0, 0)$, the positive $x$-axis extends to the right of the origin, the positive $y$-axis points upward from the origin, and the positive $z$-axis points backward (in the 3rd dimension) from the origin. Let us now provide more details of these figures.
In each figure we start from the origin, and the directed cycle is inside the unit cube defined by the two corners $(0,0,0)$ and $(1,1,1)$.

\begin{figure}[H]
\centering
   \includegraphics[width=0.3\textwidth]{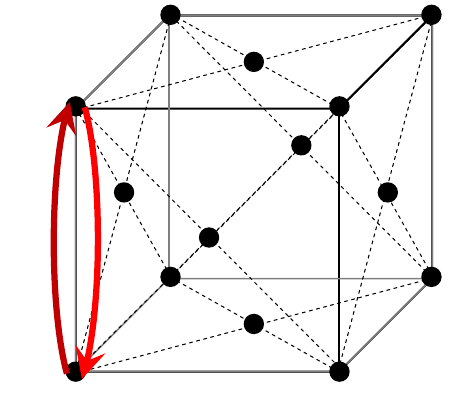} \ \includegraphics[width=0.3\textwidth]{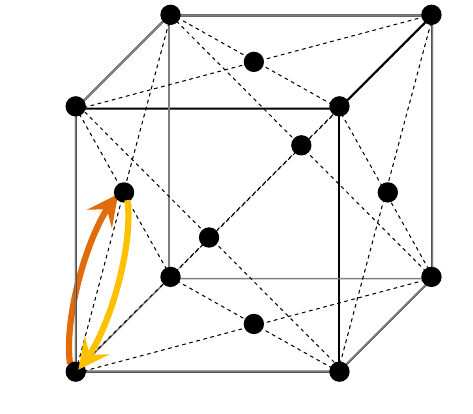}  \ \includegraphics[width=0.3\textwidth]{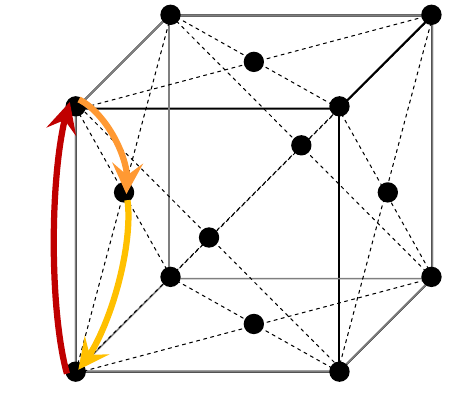}
    \caption{type 1 (left), type 2 (middle), type 3 (right).}
    \label{type1}\label{type2}\label{type3}
 \end{figure}

Figure \ref{type1} (left) tells us that $a_2+a_5=0$, where the vectors $a_2$ and $a_5$ are represented by the claret and red arrows, respectively. Here, starting from the origin and first moving along the vector $a_2$ and then moving along the vector $a_5$, we form the cycle 
 and $\{a_2, a_5\}$ is a type 1 element of the Hilbert basis which contains two vectors of the grid space. We refer to a cycle constructed from a type 1 basis element as a $1D$ base identity or as an anti-$a$ cycle.

Figure \ref{type2} (middle) shows that $d_9+d_{12}=0$, where the brown arrow represents the vector $d_9$ and the dark yellow arrow represents the vector $d_{12}$. Here starting from the origin and  first moving along the vector $d_9$ and then returning to the origin along with the vector  $d_{12}$, we form the cycle and $\{d_9,d_{12}\}$ is a type 2 element of the Hilbert basis which also contains two grid vectors. We refer to a cycle constructed from a type 2 basis element as a $2D$ base identity or as an anti-$d$ cycle.

Figure \ref{type3} (right) gives us that $a_2+d_{11}+d_{12}=0$, where the vectors $a_2, d_{11}$, and $d_{12}$ are represented by the claret, light brown and dark yellow arrows, respectively (as in the previous figures). Here, starting again from the origin we first move along the vector $a_2$, then move along the vectors  $d_{11}$ and $d_{12}$, respectively to form a cycle and, $\{a_2,d_{11},d_{12}\}$ is a type 3 element of the Hilbert basis which contains three grid vectors. We refer to a cycle constructed by a type 3 basis element as an edge-face $2D$ triangle: it represents $a-d$ dependence on 1 edge.

\begin{figure}[H]
\centering
   \includegraphics[width=0.3\textwidth]{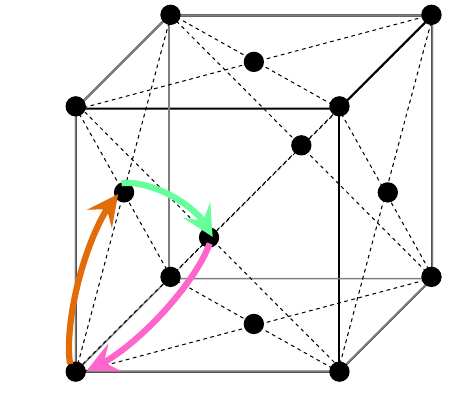} \ \includegraphics[width=0.3\textwidth]{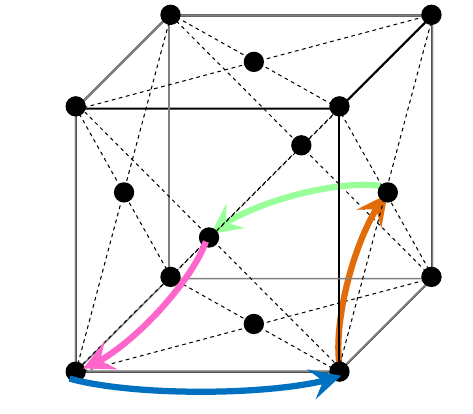}  \ \includegraphics[width=0.3\textwidth]{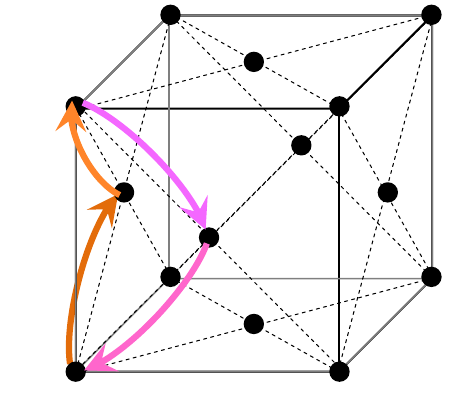}
    \caption{type 4 (left), type 5 (middle), type 6 (right).}\label{type4}\label{type5}\label{type6}
 \end{figure}

Figure \ref{type4} (left) illustrates $d_9+d_6+d_{4}=0$, where the vectors $d_9, d_6$ and $d_4$ are represented by the brown, green and pink arrows, respectively. Here, starting from the origin, we first move along the vector $d_9$ and return to the origin by moving along the vectors $d_6$ and $d_4$, respectively to form a cycle. $\{d_9,d_6,d_4\}$ is a type 4 element of the Hilbert basis which contains three grid vectors. We refer to the cycles constructed by type 4 basis elements as $3D$ triangles or corners.

Figure \ref{type5} (middle) shows that $a_1+d_9+d_8+d_{4}=0$, where the vectors $a_1, d_9, d_8$ and $d_4$ are represented by the blue, brown, light green and pink arrows, respectively. Here starting from the origin, we first move along the vector $a_1$, followed by the vectors $d_9, d_8$ and $d_4$, respectively  in an anticlockwise direction until we arrive to the starting point. $\{a_1,d_9,d_8,d_4\}$ is a type 5 element of the Hilbert basis which contains four grid vectors.

Figure \ref{type6} (right) gives us  that $d_9+d_{10}+d_2+d_4=0$, where the vectors $d_9$, and $d_{10}$ are represented by the brown, light brown arrows respectively and the vectors $d_2$, and $d_4$ are represented by light purple and pink arrows. Here starting from the origin we first move along the vector $d_9$, followed by the vectors $d_{10}, d_2$ and $d_4$, respectively, in a clockwise direction until we arrive to the origin.  $\{d_9,d_{10},d_2,d_{4}\}$ is a type 6 element of the Hilbert basis which also contains four grid vectors. We refer to a cycle constructed by a type 6 element as edge expression by 2 faces.

\begin{figure}[H]
\centering
   \includegraphics[width=0.3\textwidth]{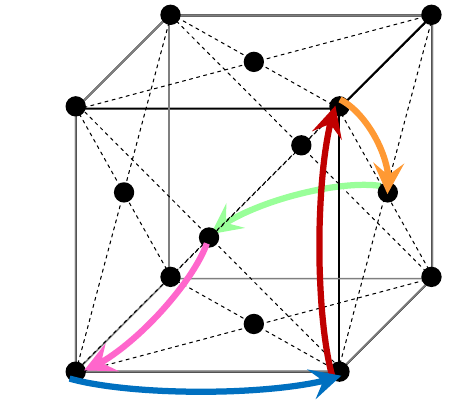} \ \includegraphics[width=0.3\textwidth]{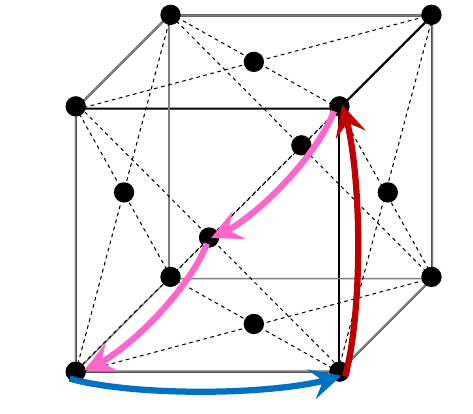}  \ \includegraphics[width=0.3\textwidth]{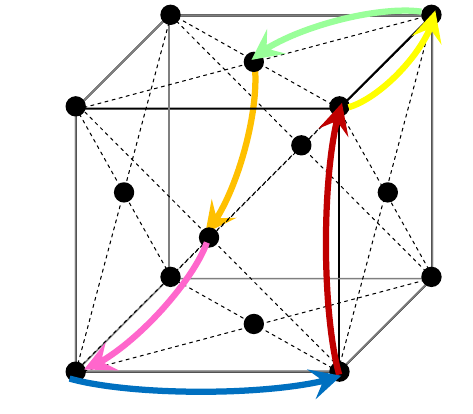}
    \caption{type 7 (left), type 8 (middle), type 9 (right).}\label{type7}\label{type8}\label{type9}
 \end{figure}

Figure \ref{type7} (left) shows that $a_1+a_{2}+d_{11}+d_8+d_4=0$, where the vectors $a_1, a_2, d_{11}, d_8$ and $d_4$ are represented by the blue, claret, light brown, light green and pink arrows, respectively. Here starting from the origin, we first move along the vector $a_1$, followed by the vectors $a_2,d_{11},d_8$, and $d_4$, respectively, 
 until we arrive to the origin. $\{a_1,a_2,d_{11},d_8,d_4\}$ is a type 7 element of the Hilbert basis which contains five grid vectors.

Figure \ref{type8}  (middle) shows that $a_1+a_2+2d_4=0$, where  the vectors $a_1,a_2$ and $d_4$ are represented by the blue, claret, and pink arrows, respectively. Here starting from origin we first move along the vector $a_1$, and then move along to vector $a_2$, followed by the vector $d_4$ twice, in that order until we arrive to the origin. This element of the Hilbert basis can be described by the vector $(1,1,0,0,0,0,0,0,0,2,0,0,0,0,0,0,0,0)^T$ in the step-space, and can be described by the support set
$\{a_1,a_2,d_4\}$  in the grid-space. It is a type 8 element of the Hilbert basis, and the directed cycle is built as an addition of 4 grid vectors. We refer to the cycles constructed by type 8 basis elements as $2D$ triangles $a-d$ dependence on 2 edges or alternatively, we may call them as face diagonals.

Figure \ref{type9}  (right) shows us that $a_1+a_2+a_3+d_{12}+d_8+d_4=0$, where the vectors $a_1,a_2,a_3,d_{12},d_8$, and $d_4$ are represented by  the blue, claret, yellow, light green, dark yellow and pink arrows, respectively. Here we first move along the vector $a_1$, followed by the vectors $a_2,a_3,d_{12},d_8$, and $d_4$, respectively,  in an anticlockwise direction until we arrive to the origin. $\{a_1,a_2,a_3,d_{12},d_8,d_4\}$ is a type 9 element of the Hilbert basis which contains six grid vectors.
We refer to the cycles constructed by type 9 basis elements as $3D$ identity of body diagonal.

\begin{figure}[H]
\centering
   \includegraphics[width=0.3\textwidth]{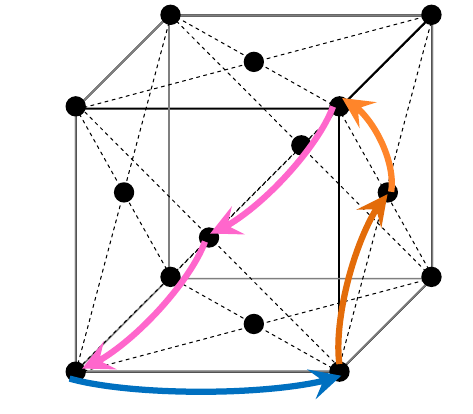} \ \includegraphics[width=0.3\textwidth]{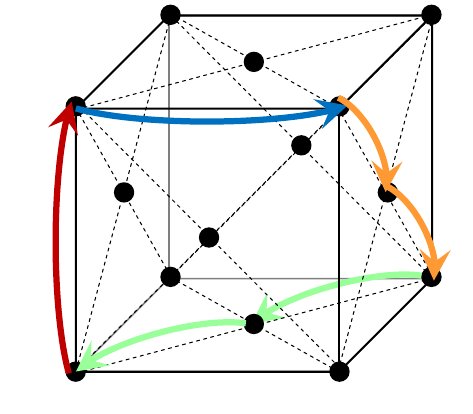}
    \caption{type 10 (left), type 11 (right).}\label{type10}\label{type11}
 \end{figure}

Figure \ref{type10} (left) tells us that $a_1+d_9+d_{10}+2d_4=0$, where $a_1,d_9,d_{10}$, and $d_4$ are represented by the blue, brown, light brown and pink colors. Here, we first move along the vector $a_1$, followed by the vectors  $d_{9},d_{10}$, and $d_4$ twice, in that order respectively, in an anticlockwise direction until we arrive to the origin. Its support set is $\{a_1,d_9,d_{10},d_4\}$ and it is a type 10 element of the Hilbert basis which is represented by a directed cycle as a sum of 5 grid vectors.

Finally, Figure \ref{type11}  (right) tells us that $a_2+a_1+2d_{11}+2d_8=0$, where the vectors $a_2,a_1,d_{11},d_8$ are represented by claret, blue, light brown and light green colors, respectively. Here we first move along to the vector $a_2$, secondly move along to the vector $a_1$, followed by the vector $d_{11}$ twice, in that order and then move the along to the vector $d_8$ twice, in that order until to reach the origin. The support set is $\{a_2,a_1,d_{11},d_8\}$ and it is a type 11 element of the Hilbert basis which means a sum of 6 vectors.

Now, we are ready to count the numbers of the possibilities in each type.

{\begin{theorem}\label{THcases}
The last column of Table \ref{tab1} shows exactly the number of different types of elements of the Hilbert basis.
\end{theorem}}
\proof{Proof goes by combinatorial arguments combined with geometric arguments for each case.

\underline{Case 1:}
A $1D$ base identity or an anti-$a$ cycle consists of two unit cube vectors such that they are opposite of each other.
We work in 3D and there are six unit cube vectors  in the matrix $A$, with half of them being opposites of each other. Therefore we have three such pairs. Alternatively, we may see these three cases geometrically: Figure \ref{type1} (left) itself illustrates one case, the other two cases could be seen if we rotate the unit cube 90 degrees around the $x$ or the $z$ axes, respectively.

\underline{Case 2:}
A $2D$ base identity or an anti-$d$ cycle consists of two diagonal vectors that are opposites of each other.  We have twelve diagonal vectors in the matrix $A$, with half of them being opposites of each other which gives six pairs.

\underline{Case 3:}
We have twelve edge-face $2D$ triangles. In Figure \ref{type1} (right) the edge-face $2D$ triangle is constructed as follows: we first choose any six of the unit cube vectors. Since each vector lies on the edge of two adjacent faces of the unit cube, there are two alternative ways to return to the origin using these faces. This validate that we have twelve triangles.   Specifically, in this figure, the unit cube vector $a_2$ is chosen first, which lies at the edge of the front and left faces. Then, the diagonal vectors $d_{11}$ and $d_{12}$ are used sequentially to return to the origin through the left face. But also after the first move, we may also use the diagonal vectors $d_2$ and $d_4$, respectively to go back to the origin through the front face.
Each edge-face $2D$ triangle can be obtained by rotating the triangle on the figure.

\underline{Case 4:}
 A $3D$ triangle or a corner has eight such cases. In the figure
 \ref{type4} (left) we may see that  by isometric transformation which brings this cube to itself such that the given corner in that figure is going to any other corner then this  $3D$ triangle goes to one of the other triangles. In fact, we have eight corners and we may have a bijection where for each corner we assign one such $3D$ triangle. One may also see that reverse ordering of the vectors in the given figure also gives us a $3D$ triangle and this triangle is assigned to the opposite  corner.

\underline{Case 5:}
There are twenty four cases of type 5 members of the Hilbert basis: From figure
 \ref{type5} (middle) we may first see that, starting from the origin, there are six possible choices to move along the unit cube vectors. After making the first move, we have four diagonal vector choices, which are perpendicular to the unit vector which we chose. Once the vectors for the first and second moves are determined, these choices uniquely define the remaining two moves. 
  Consequently, we have the only given two diagonal vectors that allow us to return to the origin.

\underline{Case 6:}
There are six cases of edge expression by 2 faces: Figure \ref{type6} (right) tells us that we have six edges in one corner and for each that edge we can order one cycle where in one face we go and in the other face we come back. We have six edges and we may have a bijection where for each edge we assign one such cycle. Note that it is clear that reverse ordering of the vectors in the given figure also gives us such cycle and we assign such cycle to the opposite edge we considered.

\underline{Case 7:}
There are twenty four cases of type 7 members of the Hilbert basis: Figure \ref{type7} (left) tells us that for the first two moves we need to choose two unit cube vectors which are perpendicular to each other. In particular, we first have six possible choices to move along the unit cube vectors and we have four choices among the remaining unit cube vectors which are perpendicular to the vector which we first used. However, we can select the same two unit vectors in the opposite order, therefore
we have $\frac{6\cdot4}{2} = 12$ choices.  At this point, on the back way, one of the diagonal vectors is fixed (the last, pink, arrow in the figure)
determined by the two chosen unit vectors. For the other two vectors we have a choice which of them is increasing and which of them is decreasing the third coordinate (that is the remaining coordinate, the one has not been changed by the chosen unit vectors), as the orange and light green arrows in the figure. Hence, all together, there are $12\cdot 2=24$ such members.

\underline{Case 8:}
There are twelve $2D$ triangles $a-d$ dependence on 2 edges. Figure \ref{type8} (middle)
tells us that we first have six possible choices to move along the unit cube vectors and
we have four choices among the remaining unit cube vectors which are perpendicular to the first chosen vector. However, we can select the same two unit vectors in the opposite order, therefore we have $\frac{6\cdot4}{2} = 12$ choices.  At this point, on the back way, the diagonal vectors are fixed (pink arrows in the figure) as their direction is determined by the two chosen unit vectors. This means we have $12$ such triangles.

\underline{Case 9:}
There are eight $3D$ identity of body diagonals. Figure \ref{type9} (right) tells us that starting from the origin, using unit cube vectors we first reach one of the eight corners: $(\pm2,\pm2,\pm2)$. When we come back we use diagonal vectors which are already determined by the unit cube vectors we used. Therefore we have eight such cases.
Each such case can be obtained by mirroring the unit cube of the figure to the unit cube having the chosen corner, respectively.

\underline{Case 10:}
There are twenty four type 10 members of the Hilbert basis: Figure \ref{type10} (left) tells us that for the first move we have six choices through unit cube vectors and once we choose one of them, since we moved on an edge, then there are four different faces to choose to go back. Once we have chosen one of these faces then we use two diagonal vectors to reach to the next corner (opposite to the origin in that face), this can be done only by that two diagonal vectors. Since our face choice already determine which face we need to use to arrive back to the origin, we have no choice left to use this diagonal vector twice. Hence we have twenty four such cases.

\underline{Case 11:}
There are twenty four type 11 members of the Hilbert basis: Figure \ref{type11} (right)
tells us that for the first two moves we need to choose one unit cube vector each which are perpendicular to each other.
 As we have already shown in some previous cases, this can be done in $12$ cases. At this point
 we arrive to a face diagonal corner. We now start going back. For this whenever we want to change the third coordinate which is the $z$ coordinate according to the figure we show, we can change it in two ways: moving along the positive $z$-axis which means using the diagonal vector $d_{11}$ twice or moving along to negative $z$-axis which means using the diagonal vector $d_{12}$ twice. Therefore, we have  $12\cdot 2$ options. Now to arrive back to the origin we have exactly one diagonal vector option left which is already defined by the previous choices. Using this diagonal vector twice have arrive back to the origin.
}
\hfill Q.E.D.\\

\section{Conclusions}
One of the most important 3D grids, the face-centered cubic grid was studied. It was shown that there are 11 different types of the 151 Hilbert basis elements. Each type is also graphically represented. Our results, apart from the theoretical, mathematical interest, can be important for various scientists working with the FCC grid including (solid state) physicists, crystallographers, etc. Our results can be applied in modelling in various phenomena in the FCC grid where optimal/shortest paths play an important role.

It is interesting that each element of the Hilbert basis in the FCC grid can be represented by adding vectors inside a unit cube, as we have displayed them on Figures \ref{type1}-\ref{type10}. In \cite{kov5}, it is shown that this property does not hold for the 4-dimensional generalization of the BCC grid. But, it remains an interesting question whether this property holds for the 3-dimensional BCC grid and other types of grids.

\end{document}